\documentclass[12pt]{article}
\usepackage{amsmath,amssymb,amsthm,latexsym}
\usepackage[dvips]{graphics}
\usepackage{epsfig, rotate}
\usepackage[english]{babel}
\usepackage{amsfonts}
\usepackage{array}
\usepackage{color}

\newtheorem{theorem}{{Theorem}}

\newtheorem{lemma}{{Lemma}}

\numberwithin{equation}{section} \numberwithin{theorem}{section}
\numberwithin{lemma}{section}\numberwithin{corollary}{section}

\usepackage%[pctex32]
{graphics}
\usepackage%[dvips]
{graphicx}

\newcommand{\be}{\begin{equation}}
\newcommand{\ee}{\end{equation}}
\newcommand{\beaa}{\begin{eqnarray*}}
\newcommand{\eeaa}{\end{eqnarray*}}
\newcommand{\bea}{\begin{eqnarray}}
\newcommand{\eea}{\end{eqnarray}}

\newcommand{\bei}{\begin{itemize}}
\newcommand{\eei}{\end{itemize}}

\newcommand{\bd}{\mathbf}

\newcommand{\td}{\overset{d}\to}
\newcommand{\tp}{\overset{p}\to}
\newcommand{\ed}{\overset{d}=}

\numberwithin{equation}{section}

\def\E{\mathbb{E}}

\def\td{\stackrel{d}{\rightarrow}}

\begin{document}

\noindent {\bf \Large  Limiting Distributions of Spectral Radii for \\
Product of Matrices from the Spherical Ensemble}

\vspace{20pt}
\noindent{\bf Shuhua Chang$^1$, Deli Li$^2$,  Yongcheng Qi$^3$}

\vspace{10pt}

{\small
\noindent $^1$Coordinated Innovation Center for Computable Modeling in Management Science, Tianjin University of Finance and Economics, Tianjin 300222, PR China. \\ Email: szhang@tjufe.edu.cn

\vspace{10pt}
\noindent $^2$Department of Mathematical Sciences, Lakehead University, Thunder Bay, ON P7B 5E1, Canada
\\ Email: dli@lakeheadu.ca

\vspace{10pt}

\noindent $^3$Department of Mathematics and Statistics, University of Minnesota Duluth,
1117 University Drive, Duluth, MN 55812, USA.
\\ Email: yqi@d.umn.edu

\date{\today}

\vspace{20pt}

\noindent{\bf Abstract.}
Consider the product of $m$ independent $n\times n$ random matrices from the spherical ensemble for $m\ge 1$.
The spectral radius is defined as the maximum absolute value of the $n$ eigenvalues of the product matrix. When $m=1$, the limiting distribution for the spectral radii has been obtained by Jiang and Qi (2017).
In this paper, we investigate the limiting distributions for the spectral radii in general. When $m$ is a fixed integer, we show that the spectral radii converge weakly to distributions of functions of independent Gamma random variables.
When $m=m_n$ tends to infinity as $n$ goes to infinity,  we show that the logarithmic spectral radii have a normal limit.

\vspace{20pt}

\noindent {\bf Keywords:}~  limiting distribution, spectral radius, spherical ensemble, product ensemble, random matrix
}

%\vspace{20pt}

\newpage

\section{Introduction}\label{intro}

In the last few decades, random matrix theory has expanded very quickly and found applications in many areas such as heavy-nuclei (Wigner, 1955), condensed matter physics (Beenakker, 1997),
%quantum mechanics (Mehta, 2004),
number theory (Mezzadri and Snaith, 2005),
wireless communications (Couillet and Debbah, 2011), and high dimensional statistics  (Johnstone (2001, 2008) and Jiang (2009)), just to mention a few. Interested readers are referred to the Oxford Handbook of Random Matrix Theory edited by Akemann, Baik and Francesco (2011) for more references and a wide range of applications in both mathematics and physics.

The study of the largest eigenvalues of Hermitian random matrices has been very active after the discovery of the so-called Tracy-Widom
distributions. For the three Hermitian matrices including Gaussian orthogonal ensemble, Gaussian unitary ensemble and Gaussian symplectic ensemble, Tracy and Widom (1994, 1996) have proved that the largest eigenvalues converge in distribution to some distributions, now known as  the Tracy-Widom laws.  Later developments in this direction can be found in Baik et al. (1999), Tracy and Widom (2002), Johansson (2007), Johnstone (2001, 2008) and Jiang (2009), and Ram\'{\i}rez et al. (2011).

The study of non-Hermitian matrices, initiated by Ginibre (1965) for Gaussian random matrices, has attracted much attention as well, and  applications are found in areas such as quantum chromodynamics, chaotic quantum systems and growth processes; see, e.g., Akemann, Baik and Francesco (2011) for more descriptions. For non-Hermitian matrices, the largest absolute values of their eigenvalues are refereed to as the  spectral radii.  Rider (2003, 2004) and Rider and
Sinclair (2014) consider the real, complex and symplectic Ginibre ensembles.
In particular, for the complex Ginibre ensemble, Rider (2003) shows that the spectral radius converges in distribution
to the Gumbel distribution. Jiang and Qi (2017) investigate the limiting distributions for the spectral radii for the spherical ensemble,
truncation of circular unitary ensemble and product of independent matrices with entries being independent complex
standard normal random variables. These limiting distributions are no longer the Tracy-Widom laws.
Gui and Qi (2018) further extend Jiang and Qi's (2017) result for the truncations of circular unitary ensemble. A common feature for all
these random matrices is the intrinsic independence structure for the absolute values of their eigenvalues, which is shared by certain determinantal point processes; see e.g.,  Hough et al. (2009).

 Let $m\ge 1$ be an integer and assume $\bd{X}_1, \cdots, \bd{X}_m$ are $m$ independent and identically distributed (i.i.d.) $n\times n$ random matrices. The product of the $m$ matrices is an $n\times n$ random matrix, denoted by
\begin{equation}\label{Xm}
\bd{X}^{(m)}=\bd{X}_1\bd{X}_2\cdots\bd{X}_m.
\end{equation}
The product of random matrices have been applied
in wireless telecommunication, disordered spin chain, the stability of large complex system, quantum transport in disordered wires, among others. See Ipsen (2015) for a survey of applications.

Some recent interests focus on the study of the limiting properties of the product ensemble $\bd{X}^{(m)}$, including the limit of the empirical spectral distributions and the spectral radii. For example, G\"{o}tze and Tikhomirov (2010),
Bordenave (2011), O'Rourke and Soshnikov (2011) and O'Rourke {\it et al.} (2015) have investigated the limiting empirical spectral distribution for the product from the complex Ginibre ensemble when $m$ is fixed, G\"otze,  K\"osters and  Tikhomirov (2015) and Zeng (2016) have obtained the limits of the empirical spectral distribution for the product from the spherical ensemble when $m$ is fixed, and Chang and Qi (2017) obtain the the limit of the empirical distributions based on scaled eigenvalues when $m=m_n$ changes with $n$.
The universality of convergence for the empirical spectral distribution is also obtained  by Bordenave (2011) and
G\"otze, K\"osters and  Tikhomirov (2015) when $m$ is a fixed integer.

When the $n^2$ entries of $\bd{X}_1$ are i.i.d. complex standard normal random variables, the limiting distribution for the spectral radii of $\bd{X}^{(m_n)}$ depends on the limits of $m_n/n$.  Three different types of limiting distributions are obtained in Jiang and Qi (2017) when $\lim_{n\to\infty}m_n/n=0$, $\lim_{n\to\infty}m_n/n=\alpha\in (0,\infty)$, and $\lim_{n\to\infty}m_n/n=\infty$.

%Eigenvalues with a joint density with a similar structure to \eqref{elastic-m} form a
%determinantal point process. See, e.g.,  Hough {\it et al.} (2009) for properties of determinantal point processes.
%Eigenvalues from the product of Ginibre ensembles and the product of truncations of independent
%Haar unitary matrices can be also modeled by determinantal point processes. By developing a special technique for determinantal point %processes, Jiang and Qi (2018) have obtained the limits for the empirical spectral distributions for the two aforementioned product %ensembles.

 Assume that $\bd{A}$ and $\bd{B}$ are two $n\times n$ random matrices and all of the $2n^2$ entries of the matrices are i.i.d. standard complex normal random variables. A spherical ensemble is defined as $\bd{X}:=\bd{A}^{-1}\bd{B}$;  see e.g., Hough {\it et al.} (2009).
Denote $\bd{z}_1, \cdots, \bd{z}_n$ as the eigenvalues of $\bd{X}$. Then  it follows from  Krishnapur (2009) that the joint probability density function of the $n$ eigenvalues is given by
\begin{equation}\label{elastic}
C_1\cdot\prod_{j<k}|z_j-z_k|^2\cdot \prod_{k=1}^n\frac{1}{(1+|z_k|^2)^{n+1}},
\end{equation}
where $C_1$ is a normalizing constant.

In this paper, we consider the product
of $m$ independent matrices from the spherical ensemble.
We are interested in the limiting distributions of the spectral radii for the product ensemble $\bd{X}^{(m)}$ when $n$ goes to infinity. We also allow that $m=m_n$ changes with $n$.

Let $\bd{X}_1, \cdots, \bd{X}_m$ be $m$ independent and identically distributed $n\times n$ random matrices that have the same distribution as $\bd{X}$ defined above.  The product ensemble $\bd{X}^{(m)}$ is defined as in \eqref{Xm}. Then we have from  Adhikari {\it et al.} (2016) that the $n$ eigenvalues $\bd{z}_1, \cdots, \bd{z}_n$ of $\bd{X}^{(m)}$ have  a joint probability density function
\begin{equation}\label{elastic-m}
C_m\cdot\prod_{j<k}|z_j-z_k|^2\cdot \prod_{k=1}^nw_m(z_k),
\end{equation}
where $C_m$ is a normalizing constant and $w_m(z)$ can be expressed in terms of a Meijer $G$-function.
%\[w_n(z)=\frac{\pi^{m-1}}{(n!)^m}G^{m,m}_{m,m}\left({{(-n,-n,\cdots, -n)_m}\atop{(0,0\cdots, 0)_m}} \Big| |z|^2\right),\]
%and $G^{m,m}_{m,m}\left({{(-n,-n,\cdots, -n)_m}\atop{(0,0\cdots, 0)_m}} \Big| |z|^2\right)$ is a Meijer $G$-function.
%See Adhikari {\it et al.}(2016).
A recursive formula for $w_m$ is obtained by Zeng (2016) as follows
\[
w_{k+1}(z)=2\pi\int^\infty_0w_k(\frac{z}{r})\frac{1}{(1+r^2)^{n+1}}\frac{d\,r}{r}
\]
for $k\ge 1$ with initial $w_1(z)=\displaystyle\frac{1}{(1+|z|^2)^{n+1}}$. Clearly, \eqref{elastic} is a special case of \eqref{elastic-m} when $m=1$.

\vspace{10pt}

When $m=1$, the limiting distribution has been obtained in Jiang and Qi (2017). In this paper, our objective is to obtain the limiting distributions for the spectral radii for the product ensemble $\bd{X}^{(m)}$ in the following two cases: (a) $m\ge 1$ is a fixed integer, and (b) $m=m_n$ tends to infinity as $n$ goes to infinity. We will show that the limiting distributions of the spectral radii can be expressed as the distributions of functions of independent Gamma random variables when $m$ is fixed, and the limiting distributions for the logarithmic spectral radii are normal when $m=m_n$ diverges as $n\to\infty$.

The rest of the paper is organized as follows.  The main results of the paper are introduced in Section~\ref{main},  and their proofs
are  given in Section~\ref{proof}.

\section{Main Results}\label{main}

We assume that the product $\bd{X}^{(m)}$ defined in \eqref{Xm} is the product
of $m$ i.i.d. random matrices from the spherical ensemble. Note that the eigenvalues $\bd{z}_1, \cdots, \bd{z}_n$ of $\bd{X}^{(m)}$
are complex random variables with the joint density distribution function given in \eqref{elastic-m}.  The spectral radius
of $\bd{X}^{(m)}$ is defined as
\begin{equation}\label{Mn}
M_n=\max_{1\le j\le n}|\bd{z}_j|.
\end{equation}

Let $\{E_{ijk}, ~i\ge 1, ~j\ge 1, ~k\ge 1\}$ be i.i.d random variables with standard exponential  distribution
(ie., Gamma($1$) distribution).  Set $\Gamma_{ij}[k_1:k_2]=\sum^{k_2}_{k=k_1}E_{ijk}$ for any $k_2\ge k_1\ge 1$, $i\ge 1$, $j\ge 1$, and denote $\Gamma_{ij}=\Gamma_{ij}[1:i]=\sum^i_{k=1}E_{ijk}$ for $i\ge 1, ~j\ge 1$.  Then
$\Gamma_{ij}, j\ge 1, i\ge 1$ are independent random variables and $\Gamma_{ij}$ has a Gamma($i$) distribution with density function $x^{i-1}e^{-x}I(x>0)/\Gamma(i)$, where $I(A)$ denotes the indicator function of set $A$, and $\Gamma(x)$ denotes the Gamma function defined as
\[
\Gamma(x)=\int^\infty_0t^{x-1}e^{-t}dt,~~x>0.
\]

We have the following two theorems on the limiting  distributions of the spectral radius $M_n$
for the product ensemble $\bd{X}^{(m)}$. The two theorems reveal two different types of limiting distributions according to whether $m$ is fixed or divergent.

\begin{theorem}\label{thm1}  Assume $m\ge 1$ is a fixed integer. Then
\begin{equation}\label{fixedm}
\frac{M_n}{n^{m/2}}\td \max_{1\le i< \infty}\frac{1}{\prod^m_{j=1}\Gamma_{ij}^{1/2}}~~\mbox{ as }n\to\infty,
\end{equation}
where $\td$ denotes convergence in distribution.
\end{theorem}

\vspace{10pt}

\begin{theorem}\label{thm2} Assume that $m=m_n\to\infty$ as $n\to\infty$.  Then we have
\begin{equation}\label{divergentm}
\frac{\log M_n-\mu_n}{\sigma_n}\td N(0,1) ~~\mbox{ as }n\to\infty,
\end{equation}
where $\mu_n=\displaystyle\frac{m_n}{2}\sum^{n-1}_{k=1}\frac1k$ and $\sigma_n^2=m_n\pi^2/24$.
\end{theorem}

\vspace{10pt}

\noindent{\bf Remark 1.} The limiting distributions are expressed in terms of functions of independent Gamma random variables in Theorem~\ref{thm1}. The random variable on the right-hand side of \eqref{fixedm} is well defined. See Lemma~\ref{existence} for a proof.

\vspace{10pt}
\noindent{\bf Remark 2.} There is no explicit form for the distribution of the random variable defined on the right-hand side of \eqref{fixedm}
except the case $m=1$.  In fact,
if we define $H_i(x)=e^{-x}\sum^{i-1}_{j=0}\frac{x^j}{j!}$ for $i\ge 1$,  then for any $x>0$
\[
P(\frac{1}{\Gamma_{i1}^{1/2}}\le x)=P(\Gamma_{i1}\ge x^{-2})=H_i(x^{-2})~~~\mbox{ for }i\ge 1,
\]
and consequently, the distribution of the random variable on the right-hand side of \eqref{fixedm} when $m=1$ is
\[
H(x)=P(\max_{1\le i< \infty}\frac{1}{\Gamma_{i1}^{1/2}}\le x)=P(\max_{1\le i<\infty}\frac{1}{\Gamma_{i1}^{1/2}}\le x)=\prod^\infty_{i=1}
H_i(x^{-2}), ~~x>0.
\]
This is exactly what Jiang and Qi (2017) have obtained in their Theorem~1. Meanwhile, they have verified that $1-H(x)\sim \frac{1}{x^2}$ as $x\to\infty$, and therefore, $H(x)$ is a heavy-tailed distribution.

\vspace{10pt}
\noindent{\bf Remark 3.}  In Theorem~\ref{thm2}, the limiting distributions are obtained for logarithmic spectral radius $\log M_n$.
It is possible to show that there do not exist real constants $a_n$ and $b_n>0$ such that $(M_n-a_n)/b_n$ converges in distribution
to a non-degenerate distribution function.

\section{Proofs}\label{proof}

First, we will introduce some notation, and then present some important lemmas. The proofs of the two main results are given afterwards.

Let $\ed$ and $\tp$ denote equality in distribution and convergence in probability. For a sequence of random variables $X_n$, $n\ge 1$
and any sequence of positive constants $a_n$, $n\ge 1$, notation $X_n=o_p(a_n)$ means $X_n/a_n\td 0$ as $n\to\infty$.
Notation $X_n=O_p(a_n)$ implies that $\lim_{c\to\infty}\limsup_{n\to\infty}P(|X_n/a_n|>c)=0$. In particular, if $X_n/a_n$ converges in distribution, then we have $X_n=O_p(a_n)$.

%Let $Y_1, \cdots, Y_n$ be $n$ independent positive random variables such that the density function of $Y_j$ is proportional to
%$y^{2j-1}w_{m_n}(y)I(y>0)$ for $1\le j\le n$, where $I(A)$ denotes the indictor function of a measurable set $A$.

Let $U_1, \cdots, U_n$ be independent random variables uniformly distributed over $(0,1)$ and define $U_{(1)}\le \cdots\le U_{(n)}$
as the order statistics of $U_1, \cdots, U_n$.

Assume that $\{s_{j, r}, ~1\le r\le m,~1\le j\le n\}$ are independent random variables, and the density of $s_{j,r}$ is proportional to $\frac{y^{j-1}}{(1+y)^{n+1}}I(y>0)$ for $1\le r\le m$, $1\le j\le n$.

Recall that $\Gamma(x)$ denotes the Gamma function.  Write $\psi(x)=\Gamma'(x)/\Gamma(x)$, $x>0$,  which is called the digamma function.
Since $\psi(x)=\frac{d}{dx}\log\Gamma(x)$, we have
\begin{equation}\label{gammatopsi}
\frac{\Gamma(b)}{\Gamma(a)}=\exp\Big(\log \Gamma(b)-\log\Gamma(a)\Big)=\exp\left(\int^b_a\psi(x)dx\right)~~\mbox{ for }a>0,~b>0.
\end{equation}

\begin{lemma}\label{mean-variance}
Let random variable  $Y$ have a Gamma($\alpha$) distribution and $X=\log(Y)$. Then the moment generating function of $X$
is given by
\[
\beta(t):=\E(e^{tX})=\frac{\Gamma(\alpha+t)}{\Gamma(\alpha)} ~~\mbox{ for } t>-\alpha.
\]
Moreover, $\E(X)=\psi(\alpha)$ and $\text{Var}(X)=\psi'(\alpha)$.
\end{lemma}

\noindent{\it Proof.}  Note $\beta(t)=\E(Y^t)$. We have for $t>-\alpha$
\[
\beta(t)=\frac{1}{\Gamma(\alpha)}\int^\infty_0y^ty^{\alpha-1}e^{-y}dy=\frac{1}{\Gamma(\alpha)}\int^\infty_0y^{\alpha+t-1}e^{-y}dy
=\frac{\Gamma(\alpha+t)}{\Gamma(\alpha)}.
\]
Then $\E(X)=\beta'(0)=\Gamma'(\alpha)/\Gamma(\alpha)=\psi(\alpha)$.
Further, we have $\E(X^2)=\beta''(0)=\Gamma''(\alpha)/\Gamma(\alpha)$. Hence we have
\[
\text{Var}(X)=\E(X^2)-\psi^2(\alpha)=\frac{\Gamma''(\alpha)}{\Gamma(\alpha)}-(\frac{\Gamma'(\alpha)}{\Gamma(\alpha)})^2
=\frac{d}{dt}\frac{\Gamma'(t)}{\Gamma(t)}\Big|_{t=\alpha}=\psi'(\alpha).
\]
This completes the proof of the lemma. \hfill$\blacksquare$

Next, we collect some properties of the bigamma function $\psi(x)$.

\begin{lemma}\label{digammafunction} For the bigamma function $\psi(x)$ we have

\noindent a. (Formulas 6.3.18 in Abramowitz and Stegun (1972))
\begin{equation}\label{psi}
\psi(x)=\log x-\frac{1}{2x}+O\Big(\frac{1}{x^2}\Big)~~\mbox{ as }~x\to \infty.
\end{equation}

\noindent b. (Formula 6.3.2 in Abramowitz and Stegun (1972))
\[
\psi(1)=-\gamma, ~~\psi(n)=-\gamma+\sum^{n-1}_{k=1}\frac1k~~~\mbox{for } n\ge 2,
\]
where $\gamma=0.57721\cdots$ is the Euler constant.

\noindent c. (Formula 6.4.10 in Abramowitz and Stegun (1972))
\[
\psi'(x)=\sum^\infty_{k=0}\frac{1}{(k+x)^2}, ~~~x>0.
\]
\end{lemma}

\vspace{10pt}

From Lemmas~\ref{digammafunction} and \ref{mean-variance} we have
\[
\psi'(1)=\sum^\infty_{k=1}\frac{1}{k^2}=\frac{\pi^2}6, ~~\psi(n)-\psi(1)=\sum^{n-1}_{k=1}\frac1k
\]
and
\[
\E(\log\Gamma_{ij})=\psi(i),~~~~\mbox{Var}(\log\Gamma_{ij})=\psi'(i).
\]
Therefore, the constants $\mu_n$ and $\sigma_n^2$ in Theorem~\ref{thm2} can be rewritten as
\[
\mu_n=m_n(\psi(n)-\psi(1))~~~\mbox{ and }~~~\sigma_n^2=m_n\pi^2/24=\frac{m_n\psi'(1)}{4}.
\]

\begin{lemma}\label{existence} For each fixed integer $m\ge 1$,
the random variable
\[
M:=\lim_{n\to\infty}\max_{1\le i\le n}\frac{1}{\prod^m_{j=1}\Gamma_{ij}^{1/2}}\]
 is well defined, and
$P(M<\infty)=1$.
\end{lemma}

\noindent{\it Proof.} Since $\max_{1\le i\le n}\frac{1}{\prod^m_{j=1}\Gamma_{ij}^{1/2}}$ is non-decreasing in $n$ with probability one,
 the limit $M$ exists and  $M>0$.  Note that
\begin{equation}\label{M4=}
M^4=\max_{1\le i< \infty}\frac{1}{\prod^m_{j=1}\Gamma_{ij}^{2}}\le \sum^\infty_{i=1}\frac{1}{\prod^m_{j=1}\Gamma_{ij}^{2}}
=\sum^\infty_{i=1}\prod^m_{j=1}\frac{1}{\Gamma_{ij}^{2}}.
\end{equation}

It follows from Lemma~\ref{digammafunction} that $\int^x_{x-2}\psi(t)dt\ge 1.5\log x$ for all large $x\ge i_0$ for some integer $i_0\ge 3$.
Therefore, it follows from Lemma~\ref{mean-variance} and equation \eqref{gammatopsi} that for $i\ge i_0$
\[
\E(\frac1{\Gamma_{ij}^2})=\frac{\Gamma(i-2)}{\Gamma(i)}
=\exp\Big(-\int^i_{i-2}\psi(t)dt\Big)
\le\exp(-1.5\log(i))
=i^{-1.5}.
\]
By using the independence of $\Gamma_{ij}$ we have
\[
\E(\sum^\infty_{i=i_0}\prod^m_{j=1}\frac{1}{\Gamma_{ij}^{2}})=\sum^\infty_{i=i_0}\E(\prod^m_{j=1}\frac{1}{\Gamma_{ij}^{2}})\le
\sum^\infty_{i=i_0}\prod^m_{j=1}\E(\frac{1}{\Gamma_{ij}^{2}})\le \sum^\infty_{i=i_0}\prod^m_{j=1}\E(\frac{1}{\Gamma_{ij}^{2}})\le
\sum^\infty_{i=i_0}i^{-1.5m}<\infty,
\]
 and hence, $P(\sum^\infty_{i=i_0}\prod^m_{j=1}\frac{1}{\Gamma_{ij}^{2}})<\infty)=1$, which together with \eqref{M4=} implies $P(M<\infty)=1$.  This completes the proof of the lemma.
\hfill$\blacksquare$

\begin{lemma}\label{zeng2016}
$g(|\bd{z}_1|^2, \cdots, |\bd{z}_n|^2)$ and
$g(\prod^m_{j=1}s_{1,j}, \cdots, \prod^m_{j=1}s_{n,j})$
have the same distribution function for any symmetric function $g(x_1,\cdots, x_n)$.
\end{lemma}

See Lemma 2.1 in Zeng (2006).

%\begin{lemma} For each $1\le i\le n$,
%\begin{equation}\label{F1}
%Y_i^2 \mbox{ and } \prod^{m}_{j=1}s_{i,j} \mbox{ are identically distributed.}
%\end{equation}
%\end{lemma}
%See Lemma 2.1 in Zeng (2016).

\begin{lemma}\label{zeng2016b} For $1\le j\le m,  1\le i\le n$,  $s_{i,j}$ and $\frac{U_{(i)}}{1-U_{(i)}}$ are identically distributed.
\end{lemma}

See the proof of Lemma 2.3 in Zeng (2016).

\begin{lemma}\label{uniform} $(U_{(1)}, U_{(2)}, \cdots, U_{(n)})$ and $(\frac{S_1}{S_{n+1}}, \frac{S_2}{S_{n+1}}, \cdots, \frac{S_n}{S_{n+1}})$
have the same joint distribution, where $S_k=\sum^k_{j=1}E_j$ for $k\ge 1$, and  $E_j,~j\ge 1$ are independent random variables with the standard exponential distribution
\end{lemma}

See, e.g., equation (2.2.1) on page 12 in Ahsanullah and Nevzorov (2015).

\begin{lemma}\label{slln} For any fixed $m\ge 1$
\[
\max_{1\le i\le n, 1\le j\le m}|\frac{\Gamma_{ij}[2:(n+1)]}{n}-1|\to 0~~~ \mbox{ with probability one}
\]
as $n\to\infty$.
\end{lemma}

\noindent{\it Proof.} For any $i, j\ge 1$, $\Gamma_{ij}[2:(n+1)]$ is the sum of $n$ i.i.d. random variables with a Gamma($1$) distribution, we have
\[
P(\max_{1\le i\le n, 1\le j\le m}|\frac{\Gamma_{ij}[2:(n+1)]}{n}-1|>\varepsilon)\le mnP(|\frac{\Gamma_{11}[2:(n+1)]}{n}-1|>\varepsilon)
\]
for any $\varepsilon>0$. Since $\Gamma_{ij}[2:(n+1)]$ is a partial sum from a sequence of i.i.d. standard exponential random variables with $\E(\Gamma_{111})=1$ and $\E(\Gamma_{111}^3)=6<\infty$, we have from Theorem 3 in Baum and Katz (1963) that
\[
\sum^\infty_{n=1}nP(|\frac{\Gamma_{11}[2:(n+1)]}{n}-1|>\varepsilon)<\infty~~~\mbox{ for any } \varepsilon>0,
\]
which implies that
\[
\sum^\infty_{n=1}P(\max_{1\le i\le n, 1\le j\le m}|\frac{\Gamma_{ij}[2:(n+1)]}{n}-1|>\varepsilon)<\infty~~\mbox{ for any } \varepsilon>0.
\]
Then the lemma follows from  Borel-Cantelli lemma.
\hfill$\blacksquare$

%\begin{lemma}\label{zeng2016-3} Set $Y_i^2=\prod^m_{j=1}s_{i, j}$ for $1\le i\le n$.
% $Y_1^2,\cdots, Y_n^2$ are stochastically increasing, that is,
%\begin{equation}\label{F3}
%P(Y_1^2\le x)\ge P(Y_2^2\le x)\ge \cdots\ge P(Y_n^2\le x)~~~x\ge 0.
%\end{equation}
%\end{lemma}
%See Lemma 2.3 in Zeng (2016).

By setting $g(x_1, \cdots, x_n)=\max_{1\le i\le n}x_i$ in Lemma~\ref{zeng2016} we have that $M_n^2=\max_{1\le i\le n}|\mathbf{z}_i|^2$
and $\max_{1\le i\le n}\prod^m_{j=1}s_{i, j}=\max_{1\le i\le n}\prod^m_{j=1}s_{n+1-i, j}$ have the same distribution.

From Lemma~\ref{uniform},  $\Gamma_{ij}/\Gamma_{ij}[1:(n+1)]=\Gamma_{ij}[1:i]/\Gamma_{ij}[1:(n+1)]$ is identically distributed as $U_{(i)}$.
Since $1-U_{(i)}$ has the same distribution as $U_{(n+1-i)}$,
we have
\[1-\frac{\Gamma_{ij}}{\Gamma_{ij}[1:(n+1)]}=\frac{\Gamma_{ij}[(i+1):(n+1)]}{\Gamma_{ij}[1:(n+1)]}\]
has the same distribution as $U_{(n+1-i)}$.  Then it follows from Lemma~\ref{zeng2016b} that
\[
\frac{\frac{\Gamma_{ij}[(i+1):(n+1)]}{\Gamma_{ij}[1:(n+1)]}}{1-(1-\frac{\Gamma_{ij}}{\Gamma_{ij}[1:(n+1)]})}=
\frac{\Gamma_{ij}[(i+1):(n+1)]}{\Gamma_{ij}}
\]
has the same distribution as $s_{n+1-i, j}$ for any $j\ge 1$. Note that $\frac{\Gamma_{ij}[(i+1):(n+1)]}{\Gamma_{ij}}$, $i\ge 1$, $j\ge 1$
are independent random variables. Therefore,
$\prod^m_{j=1}\frac{\Gamma_{ij}[(i+1):(n+1)]}{\Gamma_{ij}}$ has the same distribution as $\prod^m_{j=1}s_{n-i+1, j}$, and
$\max_{1\le i\le n}|\mathbf{z}_i|^2$ and $\max_{1\le i\le n}\prod^m_{j=1}\frac{\Gamma_{ij}[(i+1):(n+1)]}{\Gamma_{ij}}$
have the same distribution.  This implies
\begin{equation}\label{rep}
M_n\overset{d}=\max_{1\le i\le n}\prod^m_{j=1}\sqrt{\frac{\Gamma_{ij}[(i+1):(n+1)]}{\Gamma_{ij}}}.
\end{equation}
Now we define
\[
V_i=\prod^m_{j=1}\sqrt{\frac{\Gamma_{ij}[(i+1):(n+1)]}{\Gamma_{ij}}}
\]
for $1\le i\le n$.  Then we have for any $i\in\{1, \cdots, n-1\}$
\[
V_i=\prod^m_{j=1}\sqrt{\frac{\Gamma_{ij}[(i+1):(n+1)]}{\Gamma_{ij}[1:i]}}\ge \prod^m_{j=1}\sqrt{\frac{\Gamma_{ij}[(i+2):(n+1)]}{\Gamma_{ij}[1:(i+1)]}}\ed V_{i+1},
\]
which implies
\begin{equation}\label{notation}
P(V_i\le x) \mbox{ is non-decreasing in }i\in \{1, \cdots, n\}
\end{equation}
for any $x\in\mathbb{R}$.

\vspace{10pt}

\noindent{\it Proof of Theorem~\ref{thm1}.} It follows from Lemma~\ref{slln} that as $n\to\infty$
\begin{equation}\label{maxmin}
R_n:=\max_{1\le j\le m}\max_{1\le i\le n}\frac{\Gamma_{ij}[2:(n+1)]}{n}\to 1,~~r_n:=\min_{1\le j\le m}\min_{1\le i\le n} \frac{\Gamma_{ij}[2:(n+1)]}{n}\to 1
\end{equation}
with probability one.

Define for $r\ge 1$
\[
W_r=\max_{1\le i\le r}\left(\prod^m_{j=1}\sqrt{\frac{\Gamma_{ij}[(i+1):(n+1)]}{n}}
\prod^m_{j=1}\frac1{\Gamma_{ij}^{1/2}}\right).
\]
and set $Z_r=\max_{1\le i\le r}\prod^m_{j=1}\frac1{\Gamma_{ij}^{1/2}}$.
Then we have from \eqref{rep} that $M_n/n^{m/2}\overset{d}=W_n$.  To show the theorem, it suffices to prove that $W_n\to M$ with probability one.  Let $k\ge 1$ be any fixed integer. Then we have
\[
W_n\le  \max_{1\le i\le n}\prod^m_{j=1}\sqrt{\frac{\Gamma_{ij}[2:(n+1)]}{n}}Z_n\le R_n^{m/2}M,
\]
which together with \eqref{maxmin} yields
\[
\limsup_{n\to\infty}W_n\le M.
\]
For any fixed $k\ge 2$, we have for all large $n$
\[
W_n\ge W_k\ge \Big(r_n-\frac{\max_{1\le i\le k}\max_{1\le j\le m}\Gamma_{ij}[2:k]}{n}\Big)^{m/2}Z_k.
\]
Again,  in view of \eqref{maxmin} we have that $\liminf_{n\to\infty}W_n\ge Z_k$ with probability one. Hence, by letting $k\to\infty$,
and using Lemma~\ref{existence} we get that $\liminf_{n\to\infty}W_n\ge M$. Therefore, we conclude that
 $\liminf_{n\to\infty}W_n=\limsup_{n\to\infty}W_n=M$ with probability one.
\hfill$\blacksquare$

\vspace{10pt}

\noindent{\it Proof of Theorem~\ref{thm2}.}  In view of \eqref{rep} we have
\begin{equation}\label{prod}
P(\log M_n\le \mu_n+\sigma_nx)=\prod^{n}_{i=1}a_{ni}(x)
\end{equation}
for every $x\in \mathbb{R}$, where $a_{ni}(x)=P(\log V_i\le \mu_n+\sigma_nx)$ for $1\le i\le n$. Moreover,
it follows from \eqref{notation} that for each $x\in \mathbb{R}$,
\begin{equation}\label{order}
a_{ni}(x) \mbox{ is non-decreasing in } i\in\{1, \cdots, n\}.
\end{equation}

Our goal is to show that the limit on the right-hand side of \eqref{prod} is $\Phi(x)$, which is defined as the cumulative distribution
of a standard normal random variable. It suffices to show that
\begin{equation}\label{normal}
\lim_{n\to\infty}a_{n1}(x)=\Phi(x)
\end{equation}
and
\begin{equation}\label{degenerate}
\lim_{n\to\infty}\prod^n_{i=2}a_{ni}(x)=1
\end{equation}
for every $x\in \mathbb{R}$.

Note that \eqref{normal} is equivalent to
\begin{equation}\label{normal1}
\frac{\log V_1-\mu_n}{\sigma_n}\td N(0,1).
\end{equation}

For each $i\ge 1$,  $\log(\Gamma_{ij})$, $j=1, \cdots, m_n$ are i.i.d. random variables with mean $\psi(i)$ and variance $\psi'(i)$.
Then we have
\begin{equation}\label{clt12}
\frac{\sum^{m_n}_{j=1}\log(\Gamma_{ij})-m_n\psi(i)}{\sqrt{m_n\psi'(i)}}\td N(0,1)~\mbox{ as } n\to\infty
\end{equation}
by the classic central limit theorem, and as $n\to\infty$
 \begin{equation}\label{negliable}
 \frac{1}{\sqrt{m_n}}\left(\sum^{m_n}_{j=1}\log(\Gamma_{ij}[(i+1):(n+1)])-m_n\psi(n+1-i)\right)\tp 0
 \end{equation}
 since
\begin{eqnarray*}
 & &\E\left(\frac{1}{\sqrt{m_n}}(\sum^{m_n}_{j=1}\log(\Gamma_{ij}[(i+1):(n+1)])-m_n\psi(n+1-i))\right)^2\\
 &=&\frac{1}{m_n}\text{Var}(\sum^{m_n}_{j=1}\log(\Gamma_{ij}[(i+1):(n+1)]))\\
 &=&\frac1{m_n}m_n\psi'(n+1-i)\\
 &=&O(\frac1n)\\
 &\to& 0~~~\mbox{ as } n\to\infty
 \end{eqnarray*}
 from Lemma~\ref{digammafunction} (c).

Note that for $1\le i\le n$
\begin{equation}\label{sum}
\log V_i=\frac12\left(\sum^{m_n}_{j=1}\log(\Gamma_{ij}[(i+1):(n+1)])-\sum^{m_n}_{j=1}\log(\Gamma_{ij})\right).
\end{equation}

For $i=1$, we have from \eqref{clt12} and \eqref{negliable}
\begin{eqnarray*}
&&\frac{\log V_1-\mu_n}{\sigma_n}\\
&=&\frac{\sum^{m_n}_{j=1}\log(\Gamma_{1j}[2:(n+1)])-m_n\psi(n)}{\sqrt{m_n\psi'(1)}}
-\frac{\sum^{m_n}_{j=1}\log(\Gamma_{1j})-m_n\psi(1)}{\sqrt{m_n\psi'(1)}}\\
&=&-\frac{\sum^{m_n}_{j=1}\log(\Gamma_{1j})-m_n\psi(1)}{\sqrt{m_n\psi'(1)}}+o_p(1)\\
&\td & N(0,1),
\end{eqnarray*}
proving \eqref{normal1}.

For $i=2$, by using \eqref{clt12} and \eqref{negliable} and Lemma~\ref{digammafunction} (b) we get
\begin{eqnarray*}
&&\frac{\log V_2-\mu_n}{\sigma_n}\\
&=&\frac{\sum^{m_n}_{j=1}\log(\Gamma_{2j}[3:(n+1)])-m_n\psi(n-1)}{\sqrt{m_n\psi'(1)}}
-\frac{\sum^{m_n}_{j=1}\log(\Gamma_{2j})-m_n\psi(2)}{\sqrt{m_n\psi'(2)}}\sqrt{\frac{\psi'(2)}{\psi'(1)}}\\
&&-\frac{m_n(\psi(2)-\psi(1)+\psi(n-1)-\psi(n))}{\sqrt{m_n\psi'(1)}} \\
&=&-\frac{\sqrt{m_n}(1-\frac1{n-1})}{\sqrt{\psi'(1)}}+O_p(1)\\
&\tp & -\infty,
\end{eqnarray*}
which implies $1-a_{n2}(x)=P(\log V_2>\mu_n+\sigma_nx)\to 0$ as $n\to\infty$ for any $x\in\mathbb{R}$. Hence, we conclude from \eqref{order}
that $\max_{2\le i\le n}(1-a_{ni}(x))=1-a_{n2}(x)\to 0$ as $n\to\infty$.

To show \eqref{degenerate}, it suffices to show
\begin{equation}\label{limit0}
\lim_{n\to\infty}\sum^n_{i=2}(1-a_{ni}(x))=0
\end{equation}
since $1-\sum^n_{i=2}(1-a_{ni}(x))\le \prod^n_{i=1}(1-(1-a_{ni}(x)))=\prod^n_{i=2}a_{ni}(x)\le 1$.

By applying inequality $P(X>0)\le\E(e^{X})$ and noting
that all summands on the right-hand side of \eqref{sum} are independent, we have from Lemma~\ref{mean-variance} and \eqref{gammatopsi}
that
\begin{eqnarray*}
&&1-a_{ni}(x)\\
&=&P(\sum^{m_n}_{j=1}\log(\Gamma_{ij}[(i+1):(n+1)])-\sum^{m_n}_{j=1}\log(\Gamma_{ij})-m_n(\psi(n)-\psi(1))-\sqrt{m_n\psi'(1)}x>0)\\
&\le&\E\exp(\sum^{m_n}_{j=1}\log(\Gamma_{ij}[(i+1):(n+1)])-\sum^{m_n}_{j=1}\log(\Gamma_{ij})-m_n(\psi(n)-\psi(1))-\sqrt{m_n}x)\\
&=&\Big(\frac{\Gamma(n+2-i)}{\Gamma(n+1-i)}\Big)^{m_n}\Big(\frac{\Gamma(i-1)}{\Gamma(i)}\Big)^{m_n}
\exp\Big(-m_n(\psi(n)-\psi(1))-\sqrt{m_n}x\Big)\\
&=&\exp(m_n\int^1_0(\psi(n+1-i+t)-\psi(i-1+t))dt-m_n(\psi(n)-\psi(1))-\sqrt{m_n}x)\\
&=&\exp(m_n\int^1_0(\psi(n+1-i+t)-\psi(n)-(\psi(i-1+t)-\psi(1)))dt-\sqrt{m_n}x).
\end{eqnarray*}
We have $\psi(x)$ is increasing in $x>0$ since $\psi'(x)>0$ from Lemma~\ref{digammafunction}. Therefore,
$\psi(n+1-i+t)-\psi(n)\le 0$ and $\psi(i-1+t)-\psi(1)\ge \psi(i-1)-\psi(1)$ for $t\in [0,1]$ and $i\ge 3$.
Thus, we have
\[
1-a_{ni}(x)\le\exp(-m_n(\psi(i-1)-\psi(1)))=\exp(-m_n\sum^{i-2}_{k=1}\frac1k+\sqrt{m_n\psi'(1)}x)
\]
for $3\le i\le n$.  Now we choose a positive integer $i_0\ge 3$ such that $\frac12\log(i_0-1)\ge \sqrt{\psi'(1)}|x|$. Since $\sum^{i-2}_{k=1}\frac1k\ge \log(i-1)$, we get
\[
1-a_{ni}(x)\le \exp(-\frac{m_n}{2}\log (i-1))=(i-1)^{-m_n/2}, ~~~i_0\le i\le n
\]
and
\[
\sum^n_{i=i_0}(1-a_{ni}(x))\le \sum^n_{i=i_0}(i-1)^{-\frac{m_n}2}\le\int^{n-1}_{i_0-2}\frac{1}{x^{m_n/2}}dx<\frac{2}{m_n-2}(i_0-2)^{1-\frac{m_n}2}
\]
which converges to zero as $n\to\infty$. Consequently, we have as $n\to\infty$
\[
\sum^n_{i=2}(1-a_{ni}(x))\le (i_0-2)(1-a_{n2}(x))+\sum^n_{i=i_0}(1-a_{ni}(x))\to 0,
\]
which proves \eqref{limit0}. This completes the proof of the theorem. \hfill$\blacksquare$

\vspace{20pt}

\noindent{\bf Acknowledgements.} We would like to thank the reviewer whose constructive suggestions have led to improvement in the readability of the paper. Chang's research was supported in part by the Major Research Plan of the National Natural Science Foundation of China (91430108), the National Basic Research Program (2012CB955804), the National Natural Science Foundation of China (11171251), and the Major Program of Tianjin University of Finance and Economics (ZD1302). Li's research was partially supported by a grant from the Natural Sciences and Engineering Research Council of Canada (Grant \#: RGPIN-2014-05428).

\newpage

\baselineskip 12pt
\def\ref{\par\noindent\hangindent 25pt}

\end{document}